\documentclass[11pt,reqno]{amsart}
\usepackage{amsmath}
\usepackage{amssymb}

\newtheorem{thm}{Theorem}[section]

\newtheorem{defn}{Definition}[section]
\newtheorem{prop}{Proposition}[section]
\numberwithin{equation}{section}
\newtheorem{rmk}{Remark}[section]

\def\pf{{\textit {Proof:} }}

\newcommand{\mysection}[1]{\section{#1}\setcounter{equation}{0}}

\newfont{\bb}{msbm10 at 12pt}

\def\R{\hbox{\bb R}}

\def\g{{\bar g}}
\def\h{{\bar h}}
\def\k{{\lambda}}
\def\bconn{\bar{\nabla}}
\def\hconn{\breve{\nabla}}

\def\T{\mathcal P}
\def\H{\mathcal H}

\newcommand{\bal}{\begin{aligned}}      \newcommand{\eal}{\end{aligned}}
\newcommand{\ba}{\begin{array}}      \newcommand{\ea}{\end{array}}
\newcommand{\bc}{\begin{center}}     \newcommand{\ec}{\end{center}}
\newcommand{\be}{\begin{enumerate}}  \newcommand{\ee}{\end{enumerate}}
\newcommand{\beq}{\begin{eqnarray}}  \newcommand{\eeq}{\end{eqnarray}}
\newcommand{\beQ}{\begin{eqnarray*}} \newcommand{\eeQ}{\end{eqnarray*}}
\newcommand{\bi}{\begin{itemize}}    \newcommand{\ei}{\end{itemize}}
\newcommand{\bt}{\begin{tabular}}    \newcommand{\et}{\end{tabular}}
\newcommand{\bdm}{\begin{displaymath}} \newcommand{\edm}{\end{displaymath}}

\let\pa=\partial

\def\qed{\hfill{Q.E.D.}\smallskip}

\newcommand{\ls}{\setlength{\baselineskip}{12pt}
                 \setlength{\parskip}{3mm}}

\begin{document}

\title[Positive mass theorem]{Positive mass theorems for asymptotically de Sitter spacetimes}

\author{Mingxing Luo}
\address[Mingxing Luo]{Zhejiang Institute of Modern Physics, Department of Physics, Zhejiang
University, Hangzhou, Zhejiang 310027, PR China}
\email{luo@zimp.zju.edu.cn}
\author{Naqing Xie}
\address[Naqing Xie]{Institute of Mathematics, School of Mathematical Sciences, Fudan
University, Shanghai 200433, PR China}
\email{nqxie@fudan.edu.cn}
\author{Xiao Zhang}
\address[Xiao Zhang]{Institute of Mathematics, Academy of Mathematics and
System Sciences, Chinese Academy of Sciences, Beijing 100080, PR
China}
\email{xzhang@amss.ac.cn}

\date{}

\begin{abstract}
We use planar coordinates as well as hyperbolic coordinates
to separate the de Sitter spacetime into two parts. These two ways of
cutting the de Sitter give rise to two different spatial infinities. For
spacetimes which are asymptotic to either half of the de Sitter spacetime,
we are able to provide definitions of the total energy, the total
linear momentum, the total angular momentum, respectively. And we
prove two positive mass theorems, corresponding to these two sorts of
spatial infinities, for spacelike hypersurfaces whose mean curvatures are
bounded by certain constant from above.
\end{abstract}

\maketitle \pagenumbering{arabic}

\tableofcontents

\mysection{Introduction}\ls

The positive mass theorem plays a fundamental role in general
relativity. For asymptotically flat spacetimes, the positive mass
theorem was firstly proved by Schoen and Yau \cite{SY1, SY2, SY3},
and then by Witten \cite{W} using a different method. Anti-de Sitter
spacetimes are characterized by Killing spinors with imaginary
Killing constants, which have certain nice properties. Witten's
method was extended successfully to asymptotically Anti-de Sitter
spacetimes and the positive mass theorem was proved completely and
rigorously in this case \cite{Wa, CH, Z3, M, XZ}. As recent
cosmological observations indicated that our universe possibly has a
positive cosmological constant, the positive mass theorem for
asymptotically de Sitter spacetimes gains much more importance. In
spacetimes with a positive cosmological constant, a conserved
quantity was defined in \cite{AD}. Certain discussions on its
positivity and relevant issues can be found in \cite{AD, KT, S, SIT,
AMR, CGM, GM}. In particular, the conformal mass associated with the
time-like conformal Killing vector of the de Sitter spacetime was
defined in planar coordinates for asymptotically de Sitter
spacetimes and its positivity was discussed in \cite{KT}.

The de Sitter spacetime is a maximally symmetric space with positive
constant curvature. It is covered by global coordinates where each
time slice is a 3-sphere with constant curvature and has no spatial
infinity. One does not know how to define the total energy-momentum
on this slice. However, half of the de Sitter spacetime can be
covered by planar coordinates and each time slice is a
3-Euclidean space up to a conformal factor, whose second fundamental
form is proportional to the metric. Using singular coordinate
transformations, this half de Sitter spacetime can be covered by
static coordinates where cosmological horizon occurs and
time/distance switch from one to another when crossing the horizon.

It seems thus possible to build certain positive mass theorem for
half of the de Sitter spacetime. However, as pointed out by Witten
\cite{W2}, there is no positive conserved energy in de Sitter
spacetime, and the corresponding Killing vector fields to the
Lorentzian generators are timelike in some region of de Sitter
spacetime and spacelike in some other region (see also \cite{AD}).
This indicates that there should not be the positive mass theorem in
the standard sense without extra restriction. Mathematically, de
Sitter spacetimes are characterized by Killing spinors with real
Killing constants. When one tries to generalize Witten's argument to
asymptotically de Sitter spacetimes, certain essential mathematical
difficulties occur. To overcome these and to avoid using real
Killing spinors, we need to reduce an initial data set in an
asymptotically half de Sitter spacetime in planar coordinates to a
standard asymptotically flat initial data set. The total energy, the
total linear momentum and the total angular momentum can be defined
via this asymptotically flat initial data and the positive mass
theorem can be proved in this situation. The same philosophy applies
to hyperbolic coordinates which cover a different half of the de
Sitter spacetime. By using the hyperbolic positive mass theorem
proved in \cite{Z3, XZ}, we can prove the positive mass theorem in
this case. Thus we get two different half-cuts using planar
coordinates and hyperbolic coordinates and obtain two positive mass
theorems, corresponding to the asymptotically flat and
asymptotically hyperbolic positive mass theorems respectively. We
emphasize that two different half-cuts of the de Sitter spacetime
give rise to two different spatial infinities, therefore, two
different ways to measure the total energy-momentum. It will be
interesting to see how they relate to each other.

The definition of energy-momentum in planar coordinates in this
paper is slightly different from the conformal mass defined by
Kastor and Traschen \cite{KT} (see also \cite{SIT}). Let $g_0$ and
$K_0$ be the metric and the second fundamental form of the time
slice in half de Sitter spacetime in planar coordinates. For a
$\T$-asymptotically de Sitter initial data set $(M, g, K)$, $g-g_0$
on ends is used to define the total energy, the new tensor
$h=K-\sqrt{\frac{\Lambda}{3}}g$ is used to define the total linear
momentum (as well as the total angular momentum) here. But $g-g_0$
and $K-K_0$ on ends are used to define the conformal mass in
\cite{KT}. We emphasize the conformal factor $\T$ must be constant
along spacelike hypersurfaces in the definition of $\T
$-asymptotically de Sitter initial data sets. If the conformal
factor $\T$ is constant along the spacelike hypersurface and mean
curvature of the spacelike hypersurface is bounded from above by
certain constant, then the dominant energy condition for
$\T$-asymptotically de Sitter spacetimes implies that for associated
asymptotically flat spacetimes. Thus we can transfer to the
asymptotically flat case and get the positivity by using the
positive mass theorem for asymptotically flat spacetimes \cite{SY1,
SY2, SY3, W, Z1}. In \cite{KT}, Witten's method was generalized
directly to asymptotically de Sitter spacetimes in order to prove
the positivity of the conformal mass under the assumption that the
Dirac-Witten equation has solutions which are, at least through the
first correction, eigenspinors of $\gamma ^{\hat t}$ (e.g., under
(44), page 5911 \cite{KT}). In this case it does not need to use the
positive mass theorem for asymptotically flat spacetimes.

We remark that geometrical/physical properties of de Sitter
spacetime are quite different in planar coordinates and in static
coordinates. For example, de Sitter spacetime has an apparent
horizon $\{r=\k\}$ in static coordinates. Inspired by certain
physical properties in asymptotically de Sitter spacetimes in static
coordinates, some conjectured that any such spacetime with mass
greater than the mass of de Sitter has a cosmological singularity
\cite{BBM}. However, in planar coordinates, we can construct
examples which contradict to this conjecture, at least in the
``local" version. By a theorem of Corvino \cite{Co}, there exists an
asymptotically flat, scalar flat metric $\bar g_t$ in $\R ^3$ with
positive mass, which is identically Schwarzschild $\big(1+\frac{
e^{-\frac{t}{\k}}m}{2 r}\big)^4 \breve{g}$ $(m>0)$ near infinity.
Now the initial data set $(\R ^3, e^{\frac{2t}{\k}} \bar g _t, \k
^{-1} e^{\frac{2t}{\k}} \bar g _t)$ satisfies constraint equations
and provides Cauchy data for vacuum Einstein fields equations with
positive cosmological constant, which evolves into a nontrivial
vacuum spacetime identical to Schwarzschild de Sitter in planar
coordinates near spatial infinity. By the short time existence, this
spacetime is smooth and free of singularity for short time and
therefore provides a counterexample to the conjecture in this sense.

The paper is organized as follows. In Section 2, we discuss various
coordinate systems for de Sitter spacetimes. Most of them are
well-known (cf. \cite{HE}). In Section 3, we give definitions of
$\T$-asymptotically de Sitter initial data sets and the total
energy, the total linear momentum and the total angular momentum for
these initial data sets. In Section 4, we prove the positive mass
theorem for a $\T$-asymptotically de Sitter initial data set when
its mean curvature is bounded from above by certain constant. In
Section 5, we use hyperbolic coordinates to derive a positive mass
theorem for an $\H$-asymptotically de Sitter initial data set. In
Section 6, we discuss mean curvatures of time slices in certain
asymptotically de Sitter spacetimes. In Section 7, we compute the
total angular momentum for certain time slices in the Kerr-de Sitter
spacetime.

\mysection{de Sitter spacetime}\ls

The de Sitter spacetime with cosmological constant $\Lambda >0$ is
a hypersurface embedded into 5-dimensional Minkowski spacetime $\R ^{1,4}$
 \beq
-(X ^0) ^2 + (X ^1) ^2+(X ^2) ^2+(X ^3) ^2+(X ^4) ^2 =
\frac{3}{\Lambda} \label{dS}
 \eeq
with the induced metric (c.f. \cite{HE}). It is a maximally
symmetric space with constant curvature. There are several other
coordinates widely used for de Sitter spacetime. In particular, one
has the planar (inflationary) coordinates $(t, x ^i)$ and the static
coordinates $(\bar t, \bar r, \bar{\theta}, \bar{\psi})$. Denote
 \beq
\Lambda = \frac{3}{\k ^2}, \;\;\;\;\k>0.\label{C-constant}
 \eeq
The de Sitter spacetime (\ref{dS}) can be covered by the following
global coordinates:
 \beq \label{global}
 \begin{aligned}
X^0&= \k \sinh \frac{\vec{t}}{\k},\\
X ^1 &= \k \cosh \frac{\vec{t}}{\k} \sin \vec{r}\sin \vec{\theta} \cos \vec{\psi},\\
X ^2 &= \k \cosh \frac{\vec{t}}{\k} \sin \vec{r}\sin \vec{\theta} \sin \vec{\psi},\\
X ^3 &= \k \cosh \frac{\vec{t}}{\k}\sin \vec{r}\cos \vec{\theta},\\
X ^4 &= \k \cosh \frac{\vec{t}}{\k} \cos \vec{r},
 \end{aligned}
 \eeq
where $-\infty < \vec{t} < \infty$, $0\leq \vec{r} \leq \pi$, $0
\leq \vec{\theta} \leq \pi$ and $0 \leq \vec{\psi} \leq 2\pi$.
In global coordinates, the de Sitter metric is
 \beq\label{dS-gij-global}
\tilde g _{dS} = -d\vec{t} ^2 + \k^2 \cosh^2 \frac{\vec{t}}{\k}
\big(d \vec{r} ^2 +\sin ^2 \vec{r} (d \vec{\theta} ^2 +\sin ^2
\vec{\theta} d \vec{\psi} ^2) \big)
 \eeq
And $\vec{t}$-slices are 3-spheres, which do
not have spatial infinity.

The planar coordinates are defined as follows:
 \beq \label{planar+}
 \begin{aligned}
X^0&= \k \sinh \frac{t}{\k} +\frac{\k}{2}\sum _{i=1} ^3
\big(\frac{x ^i}{\k}\big)^2 e^{\frac{t}{\k}},\\
X ^i &= x^i e^{\frac{t}{\k}} \;\;\;\;(i=1,2,3),\\
X ^4 &= -\k \cosh \frac{t}{\k} +\frac{\k}{2}\sum _{i=1} ^3
\big(\frac{x ^i}{\k}\big)^2 e^{\frac{t}{\k}},
 \end{aligned}
 \eeq
where $-\infty < t < \infty$, $-\infty < x _i < \infty$. Since
 \beQ
X ^0 -X ^4 =\k e^{\frac{t}{\k}} >0,
 \eeQ
(\ref{planar+}) covers only the upper-half de Sitter spacetime. And the lower-half part is covered by
 \beq \label{planar-}
 \begin{aligned}
X^0&= -\k \sinh \frac{t}{\k} -\frac{\k}{2}\sum _{i=1} ^3
\big(\frac{x ^i}{\k}\big)^2 e^{\frac{t}{\k}},\\
X ^i &= x^i e^{\frac{t}{\k}} \;\;\;\;(i=1,2,3),\\
X ^4 &= \k \cosh \frac{t}{\k} -\frac{\k}{2}\sum _{i=1} ^3
\big(\frac{x ^i}{\k}\big)^2 e^{\frac{t}{\k}}.
 \end{aligned}
 \eeq
When $t$ goes from $-\infty$ to $\infty$,
$X ^0 -X ^4$ goes from $0$ to $\infty$ on the upper-half de Sitter spacetime
while $X ^0 -X ^4$ goes from $0$ to $-\infty$ on the lower-half part.
However, de Sitter spacetime can not be fully covered by two planar coordinates since the hypersurface
 \beQ
X ^0 - X ^4=0
 \eeQ
is excluded. In planar coordinates, the de Sitter metric is
 \beq \label{dS-gij-planar}
\tilde g _{dS} = -dt ^2 + e ^\frac{2t}{\k} \big((dx ^1)^2 +(dx ^2)^2
+(dx ^3)^2\big).
 \eeq
The second fundamental form $K$ of a $t$-slice is
 \beq \label{dS-hij-planar}
K=\frac{1}{\k} g
 \eeq
with respect to the timelike unit normal $\partial _t$.
In Minkowski spacetime $\R ^{1,4}$, $\partial _t$ is
upward in the upper-half de Sitter spacetime and downward in the lower-half part.

The static coordinates are defined as follows:
 \beq \label{static+-}
 \begin{aligned}
X ^0 &= \sqrt{\k ^2 -\bar r ^2} \sinh \frac{\bar t}{\k},\\
X ^1 &= \bar r \sin \bar \theta \cos \bar \psi,\\
X ^2 &= \bar r \sin \bar \theta \sin \bar \psi,\\
X ^3 &= \bar r \cos \bar \theta,\\
X ^4 &= -\sqrt{\k ^2 -\bar r ^2} \cosh \frac{\bar t}{\k},
 \end{aligned}
 \eeq
where $-\infty < \bar t < \infty$, $0< \bar r < \k$, $0 \leq \theta
\leq \pi$ and $0 \leq \psi \leq 2\pi$. Since
 \beQ
 \begin{aligned}
X ^0 -X ^4 &=\sqrt{\k^2 -\bar r ^2} e^{\frac{\bar t}{\k}}>0,\\
X ^0 +X ^4 &=-\sqrt{\k^2 -\bar r ^2} e^{-\frac{\bar
t}{\k}} <0,
 \end{aligned}
 \eeQ
(\ref{static+-}) covers only a quarter of the de Sitter spacetime.
The static coordinates for $\bar r >\k$ are defined as follows:
 \beq \label{static++}
 \begin{aligned}
X ^0 &= \sqrt{\bar r ^2 -\k ^2} \cosh \frac{\bar t}{\k},\\
X ^1 &= \bar r \sin \bar \theta \cos \bar \psi,\\
X ^2 &= \bar r \sin \bar \theta \sin \bar \psi,\\
X ^3 &= \bar r \cos \bar \theta,\\
X ^4 &= -\sqrt{\bar r ^2 -\k ^2} \sinh \frac{\bar t}{\k}.
 \end{aligned}
 \eeq
Since
 \beQ
 \begin{aligned}
X ^0 -X ^4 &=\sqrt{\bar r ^2 -\k^2} e^{\frac{\bar t}{\k}}
>0,\\
X ^0 +X ^4 &=\sqrt{\bar r ^2 -\k ^2} e^{-\frac{\bar t}{\k}}
>0,
 \end{aligned}
 \eeQ
(\ref{static++}) covers another quarter of the de Sitter
spacetime. (\ref{static+-}) and (\ref{static++}) cover the region
which is covered by (\ref{planar+}) excluding
 \beQ
X^0 +X ^4 =0
 \eeQ
in planar coordinates. When $\bar t$ goes from $-\infty$ to
$\infty$, $X ^0 -X ^4$ goes from $0$ to $\infty$ on the upper-half
de Sitter spacetime while $X ^0 +X ^4$ goes from $-\infty$ to
$0$ in region $\{\bar r < \k \}$ and it goes from $\infty$ to $0$ in
region $\{\bar r > \k \}$.

Static coordinates of the lower-half de Sitter $X ^0 -X ^4 <0$ can be defined in a similar way.
These four static coordinates cover the de Sitter spacetime excluding
 \beQ
X ^0 \mp X ^4=0
 \eeQ
which correspond to $\bar r=\k$ (cosmological horizons) for fixed
$\bar t$. The de Sitter metric is
 \beq \label{metric-static}
\tilde g _{dS} = -\big(1-\frac{\bar r ^2}{\k^2}\big)d\bar t ^2 +
\frac{dr ^2}{1-\frac{\bar r ^2}{\k^2}} +\bar r ^2 \big(d \bar \theta
^2 +\sin ^2 \bar \theta d \bar \psi ^2\big)
 \eeq
in the static coordinates.

The hyperbolic coordinates $(T,R,\theta,\psi)$ are defined as follows:
 \beq \label{hyperb}
 \begin{aligned}
X^0 &= \k \sinh\frac{T}{\k} \cosh \frac{R}{\k}, \\
X^1 &= \k \sinh\frac{T}{\k} \sinh \frac{R}{\k}\sin\theta\cos\psi,\\
X^2 &= \k \sinh\frac{T}{\k} \sinh \frac{R}{\k}\sin\theta\sin\psi,\\
X^3 &= \k \sinh\frac{T}{\k} \sinh \frac{R}{\k} \cos\theta,\\
X^4 &= -\k \cosh \frac{T}{\k},
 \end{aligned}
 \eeq
where $-\infty < T < \infty$, $0< R < \infty$, $0 \leq \theta \leq
\pi$ and $0 \leq \psi \leq 2\pi$. The hyperbolic coordinates cover
half the de Sitter spacetime
 \beQ
X ^4 <-\k.
 \eeQ
The de Sitter metric in the hyperbolic coordinates is
 \beq \label{metric-hyperb}
\tilde g _{dS} = -dT^2 + \sinh^2 \frac{T}{\k}\big(dR^2+\k^2\sinh^2
\frac{R}{\k}(d\theta^2+\sin^2\theta d \psi^2)\big).
 \eeq
The second fundamental form $K$ of a $T-$slice
 \beq \label{dS-hij-hyperb}
K=\frac{1}{\k}\coth \frac{T}{\k} g.
 \eeq

Now we give the relation between planar coordinates and static
coordinates. Let planar coordinates
 \beq \label{planar-pole}
 \begin{aligned}
x ^1 &= r \sin  \theta \cos \psi,\\
x ^2 &= r \sin  \theta \sin \psi,\\
x ^3 &= r \cos  \theta.
 \end{aligned}
 \eeq
By equating (\ref{planar+}) to (\ref{static+-}), and (\ref{planar+})
to (\ref{static++}), we obtain in the region $\bar r <\k$ of the upper-half de Sitter spacetime,
 \beq \label{static++planar+}
 \begin{aligned}
\bar t&= t - \frac{\k}{2} \ln \big(1-\frac{r ^2 e ^{\frac{2t}{\k}}}{\k ^2}\big),\\
\bar r&= r e ^{\frac{t}{\k}},\\
\bar \theta &=\theta,\\
\bar \psi &=\psi.
 \end{aligned}
 \eeq
In the region $\bar r >\k$ of the upper-half de Sitter spacetime, we have
 \beq \label{static+-planar+}
 \begin{aligned}
\bar t&= t - \frac{\k}{2} \ln \big(\frac{r ^2 e ^{\frac{2t}{\k}}}{\k ^2}-1\big),\\
\bar r&= r e ^{\frac{t}{\k}},\\
\bar \theta &=\theta,\\
\bar \psi &=\psi.
 \end{aligned}
 \eeq
Thus, $t$-slices in planar coordinates are not $\bar t$-slices in
static coordinates and vice versa.

\mysection{Total energy-momenta} \ls

It is not a general fact that the positive mass theorem holds true for
spacetimes with positive cosmological constant.
However, we find that, in certain special case, the mass is indeed nonnegative.
From discussions in the previous section, we realize that
cosmological horizons occur if static coordinates are used to cover the de Sitter spacetime.
Although the hypersurface
 \beQ
X^0 -X^4 =0
 \eeQ
is excluded, planar coordinates are interesting in physics. In its
early stage, the universe is generally believed to be through a
phase of inflation, which can be described by a de Sitter spacetime
in planar coordinates. On the other hand, recent cosmological
observations indicate that the expansion of the universe is
accelerating, so the future of our universe may well be again
described by a de Sitter spacetime in these coordinates. In this
section, it will be used to study the positive mass theorem for
$\T$-asymptotically de Sitter spacetimes in a mathematically
rigorous and complete way.

Suppose $(N^{1,3},\tilde{g})$ is a spacetime which satisfies the Einstein equations
with a positive cosmological constant $\Lambda$
 \beq\label{Einstein}
\tilde R _{\alpha \beta} -\frac{\tilde R}{2} \tilde g _{\alpha
\beta} +\Lambda \tilde g _{\alpha \beta} = T _{\alpha \beta}
 \eeq
where $0 \leq \alpha, \beta \leq 3$. $N ^{1,3}$ satisfies the
dominant energy condition if, for any timelike vector $W$,
 \bi
 \item[(i)] $T _{uv} W ^u W ^v \geq 0$;
 \item[(ii)] $T ^{uv} W _u $ is a non-spacelike vector.
 \ei
In a frame, we obtain
 \beq
 \bal
T ^{00} &\geq \sqrt{T ^{0i} T _i ^0}, \\
T ^{00} &\geq |T ^{\alpha \beta}|. \label{d-e}
 \eal
 \eeq

Let $(M, g, K)$ be a spacelike hypersurface with induced Riemannian metric $g$
and second fundamental form $K$.
Choose a frame $e _\alpha$ in $N ^{1,3}$ such that $e_0$ is normal to $M$ and
$e_i$ is tangential to $M$. By the Gauss and Codazzi equations,
 \beq \label{T00-T0i}
 \begin{aligned}
T_{00}&=\frac{1}{2}\big(R+(tr_g K)^2 - |K |_g ^2-2 \Lambda
\big),\\
T _{0i}&=\nabla ^j \big(K _{ij} -g _{ij} tr _g K\big),
 \end{aligned}
 \eeq
where $\nabla _i$ is the covariant derivative and $R$ is the
scalar curvature of $g$.

The tensor
 \beq
h =K -\sqrt{\frac{\Lambda}{3}} g.
\label{h}
 \eeq
will be used to define the total momenta.

Inspired by the metric and the second fundamental form
(\ref{dS-gij-planar}), (\ref{dS-hij-planar}) of a $t$-slice in de
Sitter spacetime, we can define $\T$-asymptotically de Sitter
initial data sets.
 \begin{defn}\label{defasyds}
An initial data set $(M, g, K)$ is $\T $-asymptotically de Sitter of
order $\tau >\frac{1}{2}$ if there is a compact set $M _c$ such that
$M -M _c$ is the disjoint union of a finite number of subsets $M
_1$, $\cdots$, $M _k$ - called the ``ends" of $M$ - each
diffeomorphic to $\R ^3 -B_r$ where $B _r$ is the closed ball of
radius $r$ with center at the coordinate origin. And $g$ and $h$ are
of the forms
 \beq
 \begin{aligned}
 g&= \T ^2 \bar g, \\
 h&= \T \bar h \label{g}
 \end{aligned}
 \eeq
along $M$, where $\T$ is a certain spacetime function which is a
positive constant along $M$ and tensors $\bar g $, $\bar h$ satisfy
the following conditions
 \beq \label{condition1}
\g -  \breve{g} \in C ^{2, \alpha} _{-\tau}(M),\; \h \in C ^{0,
\alpha} _{-\tau -1}(M),\;tr _{\g} (\h) \in C ^{1, \alpha} _{-\tau
-1}(M)
 \eeq
for certain $\tau >\frac{1}{2}$, $0 <\alpha < 1$, where $\breve{g}$
is the standard metric of $\R ^3$, $C ^{k,\alpha} _{-\tau}$ is
weighted H\"{o}lder spaces (c.f. \cite {B}).
 \end{defn}

Denote the metric $\breve{g} _{\T } =\T ^2 \breve{g}$. As $(M, \bar
g, \bar h)$ is asymptotically flat in the standard sense
\cite{B,Z1}, the total energy, the total linear momentum and the
total angular momentum can be defined in the standard way. Let $\{x
^i\}$ be natural coordinates of $\R ^3$, $\g _{ij}=\g(\pa _i, \pa
_j)$, $\h _{ij}=\h(\pa _i, \pa _j)$ and
 \beQ
\tilde h ^z _{ij} = \frac{1}{2}\epsilon _{i} ^{\;\;\;uv}
({\bar{\nabla}} _u \rho _z ^2 )(\h _{vj}- \g _{vj} tr _\g (\h))
 \eeQ
where $z$ is a point in $M$ and $\rho _z$ is the distance function beginning at $z$.
The tensor $\tilde h ^z _{ij}$ is trace-free and measures the rotation of the system
with respect to $z$.
It will be referred to as the local angular momentum density tensor at point $z$.
We suppose that
 \beq \label{condition2}
 \tilde h ^z  \in  C ^{0,\alpha} _{-\tau -1} (M),
 \epsilon ^{kij} {\bconn} _k (\tilde h ^z _{ij}
 -\tilde h ^z _{ji})\in
 L _{\frac{q}{2}, -\tau -2},
 \bconn ^i (\tilde h ^z _{ij}-\tilde h ^z _{ji} ) \in
 L _{\frac{q}{2}, -\tau -2}
 \eeq
for some $q>3$, where $L _{p, \tau}$ is weighted Sobolev spaces (c.f. \cite{B}).
We further assume
 \beq \label{condition3}
R \in L ^1 (M), \;T _{0i} \in L^1 (M),\;\bconn ^i \tilde h ^z _{ji}
\in L ^1 (M).
 \eeq
Under (\ref{condition1}), (\ref{condition2}) and (\ref{condition3})
the total energy $\bar E _l$, the total linear momentum $\bar P
_{lk}$ and the total angular momentum $\bar J _{lk}(z)$ with respect
to some point $z$ of the end $M _l$ are defined as follows
\cite{ADM, RegT, Z1}
 \beq
 \begin{aligned}
 \bar E _l &=\frac{1}{16\pi}\lim _{r \rightarrow \infty} \int
_{S_{r,l}}\big(\partial _j \g _{ij}-\partial _i \g _{jj}\big)*dx ^i,
\label{bEl}\\
 \bar P _{lk} &= \frac{1}{8\pi}\lim _{r \rightarrow \infty}
\int _{S_{r,l}}\big(\h_{ki}-\g _{ki} tr_{\g} (\h)\big)*dx ^i,\label{bPlk}\\
 \bar J _{lk} (z)&= \frac{1}{8\pi}\lim _{r \rightarrow \infty} \int
_{S_{r,l}}\tilde h ^z _{ki} *dx ^i \label{bJlk}.
 \end{aligned}
 \eeq
Definitions of $\bar E _l $, $\bar P _{lk}$, $\bar J _{lk} (z)$
are independent of the choice of local coordinates on ends
\cite{B,C1,Z1,Z2}.
 \begin{defn}
The total energy $E_l$, the total linear momentum $P _{lk}$ and the
total angular momentum $J _{lk}$ of the end $M _l$ are defined as
 \beq \label{definition}
E_l = \T \bar E _l,\;\;P _{lk}=\T ^2 \bar P _{lk},\;\;J _{lk}(z)=\T
^2 \bar J _{lk}(z)
 \eeq
for a $\T $-asymptotically de Sitter initial data set $(M, g, K)$.
 \end{defn}
\begin{rmk}
Certain functional spaces are used in the definition of $\T
$-asymptotically de Sitter initial data sets so that the positivity
follows from the positive mass theorem proved in \cite{Z1} for
asymptotically flat spacetimes.
\end{rmk}
\begin{rmk}
Denote $\{\breve{e} _i=\T ^{-1} \partial _i \}$ and $\{\breve{e} ^i
=\T dx ^i \}$. It is easy to check that
 \beQ
 \begin{aligned}
 E _l &=\frac{1}{16\pi}\lim _{r \rightarrow \infty} \int
_{S_{r,l}}\big[\breve{e} _j \big(g(\breve{e} _i,\breve{e}
_j)\big)-\breve{e} _i \big(g(\breve{e} _j,\breve{e} _j)
\big)\big]{*_{\breve{g} _{\T}}} {\breve{e}
^i},\label{El}\\
 P _{l} (\breve{e} _k) &= \frac{1}{8\pi}\lim _{r \rightarrow \infty}
\int _{S_{r,l}}\big[h(\breve{e} _k,\breve{e} _i) -g (\breve{e}
_k,\breve{e} _i) tr_{g} (h)\big]{*_{\breve{g} _{\T}}}
{\breve{e} ^i}.\label{Plk}\\
 \end{aligned}
 \eeQ
 \end{rmk}
 \begin{defn}
A future/past apparent horizon in a $\T $-asymptotically de Sitter
initial data set $(M, g, K)$ is a 2-sphere $\Sigma $ whose trace of
the second fundamental form $H _\Sigma $ satisfies
 \beq
\pm H _\Sigma =tr_{g _\Sigma } (h \big| _\Sigma )=tr _{g _\Sigma }
(K \big| _\Sigma ) -2 \sqrt{\frac{\Lambda }{3}}. \label{horizon}
 \eeq
 \end{defn}
 \begin{rmk}
Since it is spacelike, $\Sigma$ admits two smooth nonvanishing
outward pointing null normal vector fields $V_\pm$ with $V_+$ future
directed and $V_-$ past directed. Let $S_\pm$ be the smooth null
hypersurfaces near $\Sigma$ generated by the null geodesics with
initial tangents $V_\pm$. Then
 \beQ
\theta_\pm =div _\Sigma V _\pm = H _\Sigma  \pm tr _g (K|_\Sigma)
 \eeQ
are the null expansion of $S_\pm$ respectively which measure the
overall outward expansion of the future and past going light rays
emanating from $\Sigma$. $\Sigma$ is a physical future/past apparent
horizon if $\theta _\pm=0$ respectively. Thus the physical apparent
horizon is not the apparent horizon defined by (\ref{horizon}) (or
by (\ref{horizon-H}) in Section 5).
 \end{rmk}

In planar coordinates, the Schwarzschild-de Sitter metric can be
written as (the McVittie form, cf. \cite{KT})
 \beQ
\tilde{g}_{Sch}=-\frac{(1-\frac{m}{2Ar})^2}{(1+\frac{m}{2Ar})^2}
dt^2 +A^2 (1+\frac{m}{2Ar})^4\delta_{ij}dx^idx^j
 \eeQ
where $A(t)=e^{\frac{t}{\k}}$. It is easy to check that, for a
$t$-slice,
 \beQ
E =m, \;\;P _{k}=0,\;\;J _k (z)=0.
 \eeQ
Furthermore, if $m>0$, $\{r=\frac{m}{2A}\}$ is a minimal 2-sphere in
a $t$-slice, therefore the trace of its second fundamental form
satisfies (\ref{horizon}) and it is an apparent horizon.

\mysection{Positive mass theorem} \ls

Suppose $(M, g, K)$ is a $\T $-asymptotically de Sitter initial data
set. Let $(M, \g, p)$ be a generalized asymptotically flat initial
data set of order $\tau$, $1 \geq \tau >\frac{1}{2}$, in the sense
of \cite{Z1}, which means that $p$ is an arbitrary and not
necessarily symmetric 2-tensor. Suppose $(M, \g, p)$ satisfies the
following conditions: there exists certain compact set $\bar M _c
\supset M _c$ such that the anti-symmetric part $p^a $, $tr_{\bar
g}(p)$ are bounded on $\bar M _c$. Furthermore,
 \beQ
p \in  C ^{0, \alpha} _{-\tau -1} (M - \bar M _c),\; tr_{\bar g}(p)
\in W ^{1,\frac{q}{2}} _{-\tau -1} (M), \; \{d \theta , \;d ^*
\theta \}\in  L _{\frac{q}{2}, -\tau -2} (M)
 \eeQ
where $\theta $ is the associated 2-form of $p^a$. Suppose $(M, \bar
g, p)$ has possibly a finite number of future/past apparent
horizons $\Sigma _i$, each $\Sigma _i$ is a 2-sphere whose trace of
the second fundamental form $\bar{H} _{\Sigma _i}$ satisfies
 \beq
\pm \bar{H} _{\Sigma _i}=tr_{\g _{\Sigma _i}}(p \big| _{\Sigma _i}).
\label{horizon-p}
 \eeq
\begin{prop}
Conditions (\ref{horizon}) and (\ref{horizon-p}) are equivalent for
$p=\h$.
\end{prop}
\pf Note that $g= \T ^2 \bar g$, $h= \T \bar h$, and $\T$ is
constant along $M$ which also produces a rescalling with the same
factor on $\Sigma$. The respective mean curvatures of $\Sigma$ thus
have the relation
 \beQ
\label{2mean}H_g=\T^{-1}H_{\bar g}.
 \eeQ
On the other hand,
 \beQ
tr_g (h) = \T ^{-2} tr _{\g} (\T \h) =\T ^{-1} tr _{\g} (\h).
 \eeQ
Therefore the proposition is proved. \qed

For $(M, \g, p)$, the total energy is defined the same as $\bar E _l$
in (\ref{bEl}), the total ``linear'' momentum is defined as
 \beQ
{\bar P} _{lk} = \frac{1}{8\pi}\lim _{r \rightarrow \infty} \int
_{S_{r,l}}\big(p _{ki}-\g _{ki} tr_{\g} (p)\big)*dx ^i.
 \eeQ
However, unlike the case of symmetric 2-tensor $h$, the total
``linear'' momentum contains both translation and rotation. Denote
 \beQ
 \begin{aligned}
 \mu = \frac{1}{2}\big(\bar{R} +(tr _{\g} (p) )^2 - |p|_{\g} ^2\big),\;
 \omega _j = \bconn ^i p _{ji}-\bconn _j tr _{\g}(p),\;
 \chi _j = 2 \bconn ^i (p _{ij} - p _{ji}).
 \end{aligned}
 \eeQ
$(M, \g,p)$ satisfies the generalized dominant energy condition if
 \beq
\mu \geq \max \big\{|\omega |_{\bar g}, |\omega  + \chi |_{\bar
g}\big\}. \label{new-d-e}
 \eeq

We employ the following condition on the trace of the second fundamental form
of $M$ along $M$
 \beq
tr _g (K) \leq \sqrt{3\Lambda }. \label{trace-K}
 \eeq

\begin{prop}\label{DECs}
Under the assumption (\ref{trace-K}), the first one of the dominant
energy condition (\ref{d-e}) implies (\ref{new-d-e}) for symmetric $p=\h$.
\end{prop}
\pf  Let $\{\bar e_i\}$ be an orthonormal frame of $\bar g$ and set $K_{ij}=K(\bar e_i, \bar e_j)$.
Thus $\{e_i=\T^{-1}\bar e_i\}$ forms an orthonormal basis of $g$.
Straightforward computation yields
 \beQ
 \begin{aligned}\label{eprop2}
\mu&=\frac{1}{2}\Big(\bar R +
\T^{-2}\big(tr_{\bar g}(K)-\frac{3\T^2}{\k}\big)^2-\T^{-2}|K-\frac{\T^2}{\lambda} \bar g| _{\bar g} ^2 \Big)\\
&= \T^2 T_{00}+\frac{6\T^2}{\k ^2}-\frac{2\T^2}{\k}tr_g (K)
 \end{aligned}
 \eeQ
Under (\ref{trace-K}) and (\ref{d-e}), we have
 \beQ
 \begin{aligned}
\mu &\geq
\T^2\sqrt{\sum _{i}(\nabla_{e_j}K(e_i,e_j)-\nabla_{e_i}K(e_j,e_j))^2}\\
&=\T^{-1}\sqrt{\sum _i(\bar \nabla_{\bar e_j}K_{ij}-\bar
\nabla_{\bar e_i}K_{jj})^2}\\
&=\sqrt{\sum _i (\bar \nabla_j\bar h_{ij}-\bar \nabla_i\bar
h_{jj})^2}.
 \end{aligned}
 \eeQ \qed

Note that $\breve{g} _{ij} =\delta _{ij}$, $\breve{g} _{\T ij} =\T
^2 \delta _{ij}$. Recall the positive mass theorem proved in
\cite{Z1}.
 \begin{thm}[Zhang]\label{new-pmt-z}
Let $(M, \bar g, p)$ be a generalized asymptotically flat initial
data set of order $1 \geq \tau >\frac{1}{2}$ which has possibly a
finite number of apparent horizons. If the generalized dominant
energy condition (\ref{new-d-e}) holds, then for each end $M _l$,
 \begin{eqnarray}
\bar{E} _l \geq |\bar{P} _{l}|_{\breve{g}} . \label{new-pmt}
\end{eqnarray}
If equality holds in (\ref{new-pmt}) for some end $M _{l _0}$, then
$M$ has only one end. Furthermore, if $E _{l _0}=0$ and $\bar{g}$ is
$C ^2 $, $p$ is $C ^1$, then the following equations hold on $M$
 \beq
\overline{R} _{ijkl}+p _{ik}p _{jl}-p _{il}p _{jk}=0,\; \bconn _i p
_{jk} - \bconn _j p _{ik}=0,\; \bconn ^i (p _{ij} -p _{ji})
=0.\label{equality}
 \eeq
\end{thm}

By Theorem \ref{new-pmt-z}, or the original positive mass theorem
\cite{SY1,SY2,SY3,W, PT, Z1}, we obtain
 \begin{thm}\label{dspmtsym}
Let $(M, g, K)$ be a $\T $-asymptotically de Sitter initial data set
of order $1 \geq \tau >\frac{1}{2}$ which has possibly a finite
number of apparent horizons in spacetime $(N^{1,3},\tilde{g})$ with
positive cosmological constant $\Lambda
>0$. Suppose $N^{1,3}$ satisfies the dominant energy condition
(\ref{d-e}). If (\ref{trace-K}) holds along $M$, then for each end
$M _l$,
 \begin{eqnarray}
E _l \geq |P _{l}|_{\breve{g} _{\T}}. \label{dS-pmt}
 \end{eqnarray}
If $E_l=0$ for some end $M _{l _0}$, then
 \beq
\big(M, g, K\big) \equiv \big(\R ^3, \T ^2 \breve{g},
\sqrt{\frac{\Lambda}{3}} \T ^2 \breve{g}\big).\label{equality2}
 \eeq
Moreover, the mean curvature achieves the equality in
(\ref{trace-K}) and the spacetime $(N^{1,3}, \tilde g)$ is de Sitter
along $M$. In particular, $(N^{1,3}, \tilde g)$ is globally de
Sitter in the planar coordinates if it is globally hyperbolic.
 \end{thm}
 \pf Let $\{\bar e_i\}$ be an orthonormal frame of $\bar g$, then
$\{e_i=\T^{-1}\bar e_i\}$ forms an orthonormal basis of $g$. Since
 \beQ
|\bar{P} _{l}|_{\breve{g}} =\T ^{-1} |P _{l}|_{\breve{g} _{\T}},
 \eeQ
the positivity of mass (\ref{dS-pmt}) is a straightforward
consequence of Theorem \ref{new-pmt-z}. When $E_l=0$,  the first
equality in (\ref{equality}) implies that
 \beQ
\langle \overline{R} (\bar e_i, \bar e_j)\bar e_l, \bar e_k \rangle
_{\bar g}=-\bar h (\bar e_i, \bar e_k)\bar h (\bar e_j, \bar
e_l)+\bar h (\bar e_i, \bar e_l)\bar h (\bar e_j, \bar e_k).
 \eeQ
By the Gauss equation of $(M, g, K)$ in $(N ^{1,3},\tilde{g})$, one
has
 \beQ
 \begin{aligned}
 \tilde{R} _{ijkl}= &\langle R (e_i,e_j)e_l,e_k \rangle _g
 +K(e_i,e_k)K(e_j,e_l)-K(e_i,e_l)K(e_j,e_k)\\
=&\T^{-2}\big (-\bar h (\bar e_i, \bar e_k)\bar h (\bar e_j, \bar e_l)
+\bar h (\bar e_i, \bar e_l)\bar h (\bar e_j, \bar e_k)\big )\\
 & +K(e_i,e_k)K(e_j,e_l)-K(e_i,e_l)K(e_j,e_k)\\
=& -h _{ik} h _{jl} +h _{il} h_{jk}+\big(h _{ik}+\frac{1}{\k}g
_{ik}\big)\big(h _{jl}+\frac{1}{\k}g_{jl}\big)\\
 &-\big(h _{il}+\frac{1}{\k}g _{il}\big)\big(h
 _{jk}+\frac{1}{\k}g_{jk}\big)\\
=&\frac{1}{\k} \big(g_{ik} h _{jl} +g _{jl}h _{ik} -g _{il} h _{jk}
-g _{jk}h _{il}\big)+\frac{1}{\k ^2}\big(g _{ik}g _{jl} -g _{il} g
_{jk}\big)
 \end{aligned}
 \eeQ
along $M$. By the second equality in (\ref{equality}) (i.e., Codazzi
equations), we obtain
 \beQ
\tilde{R} _{0jkl}=0.
 \eeQ
Note that the first equality in (\ref{equality}) also implies that
$\mu =0$. From the computation in the proof of Proposition
\ref{DECs}, we have
 \beQ
T _{00} +\frac{6}{\k ^2} -\frac{2}{\k } tr _g (K) = 0.
 \eeQ
Thus the dominant energy conditions (\ref{d-e}) and (\ref{trace-K}) imply
 \beQ
 \begin{aligned}
 T _{\alpha \beta} &=0,\\
 tr_g(K)&=\frac{3}{\k}.
 \end{aligned}
 \eeQ
Therefore $N ^{1,3}$ is vacuum and
 \beq
tr _{\g} (\h)=0. \label{trace-bh}
 \eeq
The vanishing mass for $(M, \g,\h)$ under (\ref{trace-bh}) implies
that
 \beQ
\big(M, \g, \h\big) \equiv \big(\R ^3, \breve{g}, 0\big).
 \eeQ
Therefore (\ref{equality2}) holds and
 \beQ
\tilde{R} _{ijkl}=\frac{1}{\k ^2}\big(\tilde{g} _{ik}\tilde{g} _{jl}
-\tilde{g} _{il} \tilde{g} _{jk}\big)
 \eeQ
along $M$. Finally, we compute $\tilde{R} _{0j0l}$.
 \beQ
 \begin{aligned}
\tilde g ^{00} \tilde R _{0j0l} =&-\tilde g ^{ik} \tilde R _{ijkl}
+\frac{3}{\k
^2} \tilde g _{jl} \\
  =&-\frac{1}{\k ^2}\tilde{g} ^{ik} \big(\tilde{g} _{ik}\tilde{g} _{jl}
  -\tilde{g} _{il} \tilde{g} _{jk}\big)
+\frac{3}{\k ^2} \tilde{g}_{jl}\\
  =&\frac{1}{\k ^2} \tilde{g}_{jl}
 \end{aligned}
 \eeQ
along $M$. Thus $(N ^{1,3}, \tilde g)$ is de Sitter along $M$. In
particular, if it is globally hyperbolic, by the theorem of
Christodoulou and Klainerman \cite{CK}, we know that $(M, \g, \h)$
develops the Minkowski spacetime $\R^{1,3}$ which implies that $(N
^{1,3}, \tilde g)$ is de Sitter. \qed

\begin{rmk}
In general, the function $\T$ in an initial data set $(M, g, K)$ is
a nonconstant function. In this case, the positive mass theorem
(e.g. Theorem \ref{dspmtsym}) does not hold true.
\end{rmk}

We can also apply Theorem \ref{new-pmt-z} to obtain certain positive
mass theorem including the total angular momentum for asymptotically
de Sitter spacetimes.
 \begin{thm}
Let $(M, g, K)$ be a $\T $-asymptotically de Sitter initial data set
of order $1 \geq \tau >\frac{1}{2}$ which has no apparent horizon in
spacetime $(N^{1,3},\tilde{g})$ with positive cosmological constant
$\Lambda >0$. Suppose that there exists a point $z \in M$ such that
the local angular momentum tensor $\tilde h ^z $ satisfies
(\ref{condition2}), (\ref{condition3}). If (\ref{new-d-e}) holds for
$p=C_1 \h +C _2 \tilde h ^z$ where $C_1$ and $C_2$ are real constants,
then for each end $M _l$, we have
 \begin{eqnarray}
E _l \geq |C _1 P _{l} +C _2 J _{l} (z)|_{\breve{g} _{\T}}.
\label{dS-pmt2}
\end{eqnarray}
If equality holds in (\ref{dS-pmt2}) for some end $M _{l _0}$, then
$M$ has only one end. Furthermore, if $E _{l _0}=0$ and $g _{ij}$ is
$C ^2 $, $p _{ij}$ is $C ^1$, then (\ref{equality}) holds for this
$p$.
 \end{thm}

\mysection{Hyperbolic coordinates and Positive Mass Theorem} \ls

In this section, we use hyperbolic coordinates and the positive mass
theorem for asymptotically hyperbolic manifolds proved in \cite{Z3,
XZ} to derive a positive mass theorem for asymptotically half de
Sitter spacetime. Let
 \beq
\breve{g} _{\H} =dR^2+\k^2\sinh^2\frac{R}{\k}(d\theta^2+\sin^2\theta
d \psi^2). \label{gH}
 \eeq
The frame of the metric is
 \beQ
\breve{e}^\H _1=\frac{\partial}{\partial R}, \breve{e}^\H
_2=\frac{1}{\k\sinh\frac{R}{\k}}\frac{\partial}{\partial \theta},
\breve{e}^\H
_3=\frac{1}{\k\sinh\frac{R}{\k}\sin\theta}\frac{\partial}{\partial
\psi}
 \eeQ
and the coframe is
 \beQ
\breve{e} _{\H} ^1=dR, \breve{e} _{\H}
^2=\k\sinh\frac{R}{\k}d\theta, \breve{e}_{\H}
^3=\k\sinh\frac{R}{\k}\sin\theta d\psi.
 \eeQ
Denote the 4-vector $n^\nu \ (\nu=0,1,2,3)$
 \beQ
n^0=1, \ n^1=\sin\theta\cos\psi, \ n^2=\sin\theta\sin\psi, \
n^3=\cos\theta.
 \eeQ

\begin{defn}
An initial data set $(M, g, K)$ is $\H $-asymptotically de Sitter of
order $\tau >\frac{3}{2}$ if there is a compact set $M _c$ such that
$M -M _c$ is the disjoint union of a finite number of subsets $M
_1$, $\cdots$, $M _k$ - called the ``ends" of $M$ - each
diffeomorphic to $\R ^3 -B_R$ where $B _R$ is the closed ball of
radius $R$ with center at the coordinate origin. Suppose there
exists a spacetime function $T$ which is constant along $M$. Let
 \beq
h=K-\frac{\coth\frac{T}{\k}}{\k}g. \label{h-H}
 \eeq
And $g$ and $h$ are of the forms
 \beq
 \begin{aligned}
g&=\H^2 \bar g,\\
h&=\H \bar h \label{g-H}
 \end{aligned}
 \eeq
along $M$, where $\H=\sinh\frac{T}{\k}$, and
 \beQ
\bar g(\breve{e}^\H _i,\breve{e}^\H _j)-\breve{g}(\breve{e}^\H
_i,\breve{e}^\H _j)= a_{ij}, \ \bar h(\breve{e}^\H _i,\breve{e}^\H
_j)=\bar h_{ij}
 \eeQ
have the following asymptotic behavior on the end:
 \beq
 \begin{aligned}
a_{ij}=O(e^{-\frac{\tau}{\k}R}), \ \hconn ^\H _k
a_{ij}=O(&e^{-\frac{\tau}{\k}R}), \ \hconn ^\H _l \hconn ^\H _k
a_{ij}=O(e^{-\frac{\tau}{\k}R}),\\
\bar h _{ij}=O(e^{-\frac{\tau}{\k}R}), \ &\hconn ^\H _k \bar h
_{ij}=O(e^{-\frac{\tau}{\k}R}),
 \end{aligned}
 \eeq
where $\hconn ^\H$ is the Levi-Civita connection of the hyperbolic
metric $\breve{g}_\H$. Moreover,
 \beq
(\bar R +\frac{6}{\k^2})e^{\rho_z} \in L^1(M), \ \big(\bar \nabla^j
\bar h_{ij}-\bar \nabla_i tr_{\bar g}(\bar h)\big)e^{\rho_z} \in
L^1(M)
 \eeq
for some $z \in M$. Here $\bar R$, $\bar \nabla$, $\rho_z$ are
scalar curvature, Levi-Civita connection and distance function of
$\bar g$ respectively.
 \end{defn}
 \begin{defn}
A future/past apparent horizon in an $\H $-asymptotically de Sitter
initial data set $(M, g, K)$ is a 2-sphere $\Sigma $ whose trace of
the second fundamental form $H _\Sigma $ satisfies
 \beq \label{horizon-H}
\pm H _\Sigma =tr _{g _\Sigma } (K \big| _\Sigma )
-2\sqrt{\frac{\Lambda }{3}} \tanh\frac{T}{2\lambda}.
 \eeq
 \end{defn}
Let
 \beq
\tilde{h}=\bar h+\frac{1}{\k} \bar g.
 \eeq
The condition (\ref{horizon-H}) implies that
 \beq
\pm \bar{H} _\Sigma =tr _{\g _\Sigma } (\h \big| _\Sigma )
+2\sqrt{\frac{\Lambda }{3}}=tr_{\g _\Sigma } (\tilde{h} \big|
_\Sigma )
 \eeq
which shows that $\Sigma$ is an apparent horizon of $(M,\g, \h)$
\cite{Z3, XZ}. Denote
 \beq
\mathcal{E} _l =\hconn ^{\H,j} \bar g_{1j}- \hconn ^{\H} _1 tr
_{\breve{ g} _\H}(\bar g) +\frac{1}{\k}(a_{22}+a_{33})+2(\bar h
_{22}+\bar h _{33}).
 \eeq
The total energy-momentum of the end $M _l$ are defined as
 \beq
E^{\H} _{l \nu}=\frac{\H ^2}{16\pi}\lim_{R \rightarrow
\infty}\int_{S_{R,l}}\mathcal{E} _l n^\nu e^{\frac{R}{\k}} \breve{e}
_{\H} ^2\wedge \breve{e} _{\H} ^3
 \eeq
where $0\leq \nu \leq 3$, $S_{R,l}$ is the coordinate sphere of
radius $R$ in the end $M_l$.

Now we study the relation between the dominant energy condition of
the $\H$-asymptotically de Sitter initial data set and its
associated asymptotically hyperbolic initial data set. Suppose that
$(M,\bar g,\bar h)$ is an asymptotically hyperbolic initial data set
of order $\tau > \frac{3}{2}$. Denote
 \beQ
 \mu = \frac{1}{2}\big(\bar{R} +(tr _{\g} (\tilde{h}) )^2 - |\tilde{h}|_{\g}
 ^2\big),\
 \omega _j = \bconn ^i \tilde{h} _{ji}-\bconn _j tr _{\g}(\tilde{h}).
 \eeQ
Since
 \beQ
 \begin{aligned}
 tr _{\g} (\tilde{h})=&tr _{\g} (\bar{h})+\frac{3}{\lambda}\\
                     =&\H tr _{g} (h) +\frac{3}{\lambda}\\
                     =&\H tr _{g} (K) -\frac{3}{\lambda} \cosh{\frac{T}{\lambda}}+\frac{3}{\lambda},\\
|\tilde{h}|_{\g} ^2 =&(\bar h _{ij} +\frac{1}{\lambda} \bar g _{ij})
                      (\bar h _{kl} +\frac{1}{\lambda} \bar g _{kl})\bar g ^{ik} \bar g ^{jl}\\
                    =&\H ^2 |h|_{g} ^2 +\frac{2\H}{\lambda} tr _g (h) +\frac{3}{\lambda ^2}\\
                    =&\H^2 \big(|K| _g ^2 -\frac{2}{\lambda} \coth{\frac{T}{\lambda}} tr _g (K)+\frac{3}{\lambda ^2}\coth ^2 \frac{T}{\lambda}\big)\\
                    &+\frac{2\H}{\lambda} tr _g (K)-\frac{6}{\lambda ^2}\cosh\frac{T}{\lambda} +\frac{3}{\lambda ^2}\\
                    =&\H^2 |K| _g ^2 +\frac{2}{\lambda}\sinh\frac{T}{\lambda}(1-\cosh\frac{T}{\lambda})tr _g (K)\\
                    &+\frac{3}{\lambda ^2}\cosh^2 \frac{T}{\lambda}-\frac{6}{\lambda ^2}\cosh \frac{T}{\lambda}+\frac{3}{\lambda ^2},
 \end{aligned}
 \eeQ
we obtain
 \beQ
 \begin{aligned}
\mu =&\H ^2 T _{00} +\frac{2}{\lambda} tr_g(K) \sinh\frac{T}{\lambda} \big(1-\cosh\frac{T}{\lambda}\big)\\
     &+\frac{3}{\lambda^2}\sinh ^2 \frac{T}{\lambda} +\frac{3}{\lambda^2}\cosh ^2 \frac{T}{\lambda}-
\frac{6}{\lambda^2}\cosh \frac{T}{\lambda}+\frac{3}{\lambda ^2}\\
    =&\H ^2 T _{00} +\frac{2}{\lambda} \big(tr_g(K) \sinh\frac{T}{\lambda}-\frac{3}{\lambda} \cosh\frac{T}{\lambda} \big)\big(1-\cosh\frac{T}{\lambda}\big),\\
\omega _j =&\H ^2 T _{0i}.
 \end{aligned}
 \eeQ
Therefore, if
 \beq
tr _g (K) \sinh\frac{T}{\lambda} \leq \sqrt{3\Lambda } \cosh\frac{T}{\lambda}, \label{traceH-K}
 \eeq
then the dominant energy condition (\ref{d-e}) implies that
 \beq
\mu \geq |\omega |_{\g}. \label{h-doe}
 \eeq

\begin{thm}\label{dspmtsymhy}
Let $(M, g, K)$ be an $\H $-asymptotically de Sitter initial data
set of order $\tau >\frac{3}{2}$ which has possibly a finite number
of apparent horizons in spacetime $(N^{1,3},\tilde{g})$ with
positive cosmological constant $\Lambda
>0$. Suppose $N^{1,3}$ satisfies the dominant energy condition
(\ref{d-e}). If (\ref{traceH-K}) holds along $M$, then for each end $M
_l$,
 \beq
E^\H _{l0}\geq \sqrt{(E^\H _{l1})^2 +(E^\H _{l2} )^2 +(E^\H _{l3}) ^2}.\label{ds-pmt3}
 \eeq
If $E^\H _{l0}=0$ for some end $M _{l _0}$, then
 \beq
\big(M, g, K\big) \equiv \big(\hbox{\bb H} ^3, \sinh^2 \frac{T}{\lambda} \breve{g}_{\H},
\sqrt{\frac{\Lambda}{3}} \sinh\frac{T}{\lambda}\cosh\frac{T}{\lambda}\breve{g}_{\H} \big).\label{equality3}
 \eeq
Moreover, the mean curvature achieves the equality in
(\ref{traceH-K}) and the spacetime $(N^{1,3}, \tilde g)$ is de
Sitter along $M$. In particular, $(N^{1,3}, \tilde g)$ is globally
de Sitter in the hyperbolic coordinates if it is globally
hyperbolic.
\end{thm}
 \pf The first part of the theorem, i.e., inequality (\ref{ds-pmt3}),
is straightforward since the condition (\ref{h-doe}) ensures the
hyperbolic positive mass theorem proved in \cite{Z3,XZ}. If $E^\H
_{l0}=0$ for some end $M _{l _0}$, then $\mu =0$. This implies that
the mean curvature achieves the equality in (\ref{traceH-K}). Thus,
 \beQ
 \begin{aligned}
tr _{\g} (\tilde{h}) =&\frac{3}{\lambda},\\
|\tilde{h}|_{\g} ^2 =&\H ^2 |h|_{g} ^2 +\frac{3}{\lambda ^2}.
 \end{aligned}
 \eeQ
Therefore
 \beQ
\bar R =-tr _{\g} (\tilde{h}) ^2 +|\tilde{h}|_{\g} ^2 \geq
-\frac{6}{\lambda ^2}.
 \eeQ
Under this condition, that $E^\H _{l0}=0$ gives that
 \beQ
\big(M, \g\big) \equiv \big(\hbox{\bb H} ^3, \breve{g}_{\H}\big).
 \eeQ
Thus
 \beQ
\bar R=-\frac{6}{\k^2}+\H^2|h|^2_g= -\frac{6}{\k^2}.
 \eeQ
We obtain
 \beQ
 \begin{aligned}
 h&=0,\\
 K&=\frac{\cosh\frac{T}{\k}}{\k\sinh\frac{T}{\k}}g.
 \end{aligned}
 \eeQ
By the Gauss equation,
 \beQ
 \begin{aligned}
 \tilde{R}_{ijkl}&=R_{ijkl}+K_{ik}K_{jl}-K_{il}K_{jk}\\
                 &=-\frac{1}{\k^2
                 \sinh^2\frac{T}{\k}}\big(\tilde{g}_{ik}\tilde{g}_{jl}-\tilde{g}_{il}\tilde{g}_{jk}\big)
                 +\frac{\cosh^2
                 \frac{T}{\k}}{\k^2\sinh^2\frac{T}{\k}}\big(\tilde{g}_{ik}\tilde{g}_{jl}-\tilde{g}_{il}\tilde{g}_{jk}\big)\\
                 &=\frac{1}{\k^2}\big(\tilde{g}_{ik}\tilde{g}_{jl}-\tilde{g}_{il}\tilde{g}_{jk}\big).
 \end{aligned}
 \eeQ
By the Codazzi equations,
 \beQ
\tilde{R}_{0jkl}=0.
 \eeQ
The dominant energy condition and (\ref{traceH-K}) also imply that
 \beQ
T _{\alpha \beta}=0.
 \eeQ
So
  \beQ
 \begin{aligned}
\tilde g ^{00} \tilde R _{0j0l} =&-\tilde g ^{ik} \tilde R _{ijkl}
+\frac{3}{\k
^2} \tilde g _{jl} \\
  =&-\frac{1}{\k ^2}\tilde{g} ^{ik} \big(\tilde{g} _{ik}\tilde{g} _{jl}
  -\tilde{g} _{il} \tilde{g} _{jk}\big)
+\frac{3}{\k ^2} \tilde{g}_{jl}\\
  =&\frac{1}{\k ^2} \tilde{g}_{jl}
 \end{aligned}
 \eeQ
along $M$. Thus $(N ^{1,3}, \tilde g)$ is de Sitter along $M$ and we
prove the second part of the theorem. \qed

\begin{rmk}
In general, the function $\H$ in an initial data set $(M, g, K)$ is
a nonconstant function. In this case, the positive mass theorem
(e.g. Theorem \ref{dspmtsymhy}) does not hold true.
\end{rmk}

\begin{rmk}
We can also discuss the total angular momentum along the line of
\cite{Z3} in this case.
\end{rmk}

\mysection{Mean curvatures of hypersurfaces} \ls

In this section, we discuss mean curvatures of $\T$-asymptotically
de Sitter spacelike hypersurfaces. In particular, we discuss the
existence of constant mean curvature spacelike hypersurfaces. This
is analogous to the existence of maximal spacelike hypersurfaces in
asymptotically flat spacetimes, which was studied extensively by
Bartnik, etc (c.f. \cite{B84, BCM} and references therein). For
consideration of length of the paper, we study only the simple case
that the spacetime is globally hyperbolic whose existence is implied
by Christodoulou and Klainerman's theorem \cite{CK}. We will address
elsewhere by generalizing the methods in \cite{B84, BCM} to study
the existences of constant mean curvature in general spacetimes.

Let $\tilde g$ be the metric of an asymptotically flat spacetime
defined on $\R \times \R ^3$
 \beq
 \tilde g = -a ^2 (t,x)dt ^2 +g _{ij} (t,x) dx ^i dx ^j  \label{g1}
 \eeq
with $a>0$. The $t$-slice has the metric $g$ and the second
fundamental form
 \beQ
k _{ij} =\frac{1}{2a} \partial _t g _{ij}.
 \eeQ
Suppose $F(t,x)$ is a smooth function on $\R \times \R ^3$. Consider
a new metric
 \beq
 \tilde g _{\lambda} = -a ^2 (t,x)dt ^2 +e ^{2 F} g _{ij} (t,x) dx ^i dx ^j  \label{g2}
 \eeq
The $t$-slice has the metric $g _{\lambda} =e ^{2F} g$ and the
second fundamental form
 \beQ
K _{ij} =e^{2F} k _{ij} +\frac{1}{a} e ^{2F} (\partial _t F) g _{ij}.
 \eeQ
If we take
 \beQ
F(t,x) = \frac{1}{3} \int _0 ^t \big(\Theta -tr_g (k)\big) a (s,
x)ds
 \eeQ
with certain prescribed function $\Theta (t,x)$, then the mean
curvature of the $t$-slice in (\ref{g2}) is
 \beq
tr _{g_\lambda} (K)= \Theta (t,x). \label{mean-curvature2}
 \eeq

In their celebrated work \cite{CK}, Christodoulou and Klainerman
proved the global existence of an asymptotically flat metric for
(\ref{g1}), which is vacuum and foliated by maximal spacelike
hypersurfaces, i.e., $tr_g (k)=0$. Taking
 \beQ
\Theta = \sqrt{3 \Lambda}
 \eeQ
in (\ref{mean-curvature2}), we obtain an existence of constant mean
curvature spacelike hypersurfaces. Using estimates on the
lapse function \cite{CK}, we can re-write the metric (\ref{g2}) as
follows
 \beq
 \tilde g _{\lambda} = -a ^2 (t,x)dt ^2 +e ^{\frac{2t}{\k}}
 \hat{g} _{ij} (t,x) dx ^i dx ^j . \label{g3}
 \eeq
Therefore, the metric $\tilde g _{\lambda}$ is asymptotically de
Sitter. However, $\tilde g _{\lambda}$ does not satisfy the vacuum
Einstein fields equations with positive cosmological constant in
general. It will be interesting to extend Christodoulou and
Klainerman's work to the current case.

Finally, we note that suitable choices of $\Theta $ will give rise to
spacelike hypersurfaces with mean curvature violating condition (\ref{trace-K}).

\mysection{Kerr-de Sitter}

In this section, we compute the total angular momentum for
suitable time slices in the Kerr-de Sitter spacetime. In the
Boyer-Lindquist coordinates $(\bar{t}, \bar{r}, \bar{\theta},
\bar{\psi})$, the Kerr-de Sitter metric is
 \beQ
 \begin{aligned}
\tilde{g}_{KdS}=&-\frac{\Delta_{\bar{r}}}{U}\big(d\bar{t}-\frac{a}{\xi}\sin^2\bar{\theta}
d\bar{\psi}\big)^2+\frac{U}{\Delta_{\bar{r}}}d\bar{r}^2
+\frac{U}{\Delta_{\bar\theta}}d\bar{\theta}^2\\
&+\frac{\Delta_{\bar \theta} \sin ^2\bar{\theta}}{U}\big(a d
\bar{t}-\frac{(\bar{r}^2+a^2)}{\xi}d\bar{\psi}\big)^2,
 \end{aligned}
 \eeQ
where
 \beQ
 \begin{aligned}
\Delta_{\bar{r}}&=(\bar{r}^2+a^2)(1-\frac{\bar{r}^2}{\lambda^2})-2m\bar{r},\\
\Delta_{\bar \theta}&=1+\frac{a^2\cos^2\bar{\theta}}{\lambda^2},\\
U&=\bar{r}^2+a^2\cos^2\bar{\theta}, \\
\xi&=1+\frac{a^2}{\lambda^2}.
 \end{aligned}
 \eeQ

In Boyer-Lindquist coordinates, $m=0$ does not imply directly that
the metric is de Sitter. Thus, inspired by \cite{HT}, we employ the
following coordinate transformation
 \beQ
 \begin{aligned}
\hat{t}&=\bar{t},\\
\hat{r}\cos\hat{\theta}&=\bar{r}\cos\bar{\theta},\\
(1+\frac{a^2}{\lambda^2})\hat{r}^2&=\bar{r}^2+a^2\sin^2\bar{\theta}+\frac{a^2}{\lambda^2}\bar{r}^2\cos^2\bar{\theta},\\
\hat{\psi}&=(1+\frac{a^2}{\lambda^2})\bar{\psi}-\frac{a}{\lambda^2}\bar{t}.
 \end{aligned}
 \eeQ
In coordinates $(\hat{t}, \hat{r}, \hat{\theta}, \hat{\psi})$, the
Kerr-de Sitter metric can be written as
 \beQ
\tilde{g}_{KdS}=-\big(1-\frac{\hat{r} ^2}{\lambda ^2}\big)d\hat{t}
^2 +\frac{d\hat{r}^2}{1-\frac{\hat{r} ^2}{\lambda ^2}}+\hat{r}^2
(d\hat{\theta}^2 +\sin^2 \hat{\theta} d \hat{\psi} ^2)
+a_{\bar{\mu}\bar{\nu}}d\bar{x}^{\mu}d\bar{x}^{\nu}
 \eeQ
where the nonzero $a_{\bar{\mu}\bar{\nu}}$ are
 \beQ
 \begin{aligned}
a_{\bar{t}\bar{t}}&=\frac{2m \bar r}{U},\\
a_{\bar{t}\bar{\psi}}&=-\frac{2m a \bar{r}\sin^2\bar{\theta}}{U},\\
a_{\bar{r}\bar{r}}&=\frac{2m\bar{r}U}{(\Delta_{\bar{r}} +2m \bar r)\Delta_{\bar{r}}},\\
a_{\bar{\psi}\bar{\psi}}&=\frac{2m a^2 \bar{r} \sin ^4
\bar{\theta}}{U}.
 \end{aligned}
 \eeQ
The new coordinates $(\hat{t}, \hat{r}, \hat{\theta}, \hat{\psi})$
is indeed the static coordinates for the Kerr-de Sitter metric. Now we
transfer it into the planar coordinates. Let $(t, r, \theta, \psi)$
be polar coordinates corresponding to the planar coordinates.
The transformations are given as follows:
 \beQ
 \begin{aligned}
\hat{t} &=t -\frac{\k}{2} \ln \big|1-\frac{r ^2 A ^2}{\k ^2}\big|,\\
\hat{r}&=Ar,\\
\hat{\theta}&=\theta,\\
\hat{\psi}&=\psi,
 \end{aligned}
 \eeQ
where $A=e^{\frac{t}{\lambda}}$. In polar coordinates, the Kerr-de Sitter is
 \beQ
\tilde{g}_{KdS}=-dt^2 + e^{\frac{2t}{\lambda}}\big(dr^2
+r^2(d\theta^2 +\sin^2 \theta d\psi ^2)\big)+a_{\mu\nu}dx^{\mu}
dx^{\nu}
 \eeQ
and the nonzero components of $a_{\mu \nu}$ have the following
asymptotic behaviors£º
 \beQ
 \begin{aligned}
a_{tt}&=\frac{2m\lambda^2}{r^3A^3}B^{-\frac{3}{2}}+O(r^{-4}),\\
a_{t r}&=\frac{2m\lambda^3}{r^4A^3}B^{-\frac{3}{2}}+O(r^{-5}),\\
a_{t \theta}&=\frac{2m a^2 \lambda }{r^3A^4}B^{-\frac{5}{2}}\sin\theta\cos\theta + O(r^{-4}),\\
a_{rr}&=\frac{2m\lambda^2}{r^3A}B^{-\frac{5}{2}}+O(r^{-4}),\\
a_{r \theta}&=\frac{2m a^2 \lambda^2}{r^4A^3}\sin\theta\cos\theta B^{-\frac{5}{2}}+O(r^{-5}), \\
a_{r\psi}&=\frac{2m\lambda a \sin^2\theta}{r^2A}B^{-\frac{5}{2}}+O(r^{-3}),\\
a_{\theta \theta}&=\frac{2m a^4\sin^2\theta\cos^2\theta}{r^3A^3}B^{-\frac{5}{2}}+O(r^{-4}),\\
a_{\psi\psi}&=\frac{2m a^2
\sin^4\theta}{rA}B^{-\frac{5}{2}}+O(r^{-2})
 \end{aligned}
 \eeQ
where $B=1+\frac{a^2}{\lambda^2}\sin^2\theta$. The second
fundamental form of the $t$-slice can be computed in terms of the
formula
 \beQ
K _{ij} =\frac{1}{2 N} \big(\nabla _i N _j +\nabla _j N _i -\partial
_t \tilde g _{ij} \big)
 \eeQ
and the tensor $\bar h$ has the following asymptotic behaviors
 \beQ
 \begin{aligned}
\h_{rr}&=\frac{2 m \lambda ^2 -m a^2 \sin^2 \theta}{A ^2 B
^{\frac{5}{2}} \lambda r^3}+O(r^{-4}),\\
\h_{r \theta}&=O(r^{-4}),\\
\h_{r \psi }&=\frac{3m a \sin ^2\theta}{A^2 B^{\frac{5}{2}} r^2}+O(r^{-4}),\\
\h_{\theta \theta}&=-\frac{m \lambda}{A ^2 B ^{\frac{3}{2}} r}+O(r^{-3}),\\
\h_{\theta \psi}&=O(r^{-6}),\\
\h_{\psi \psi}&=\frac{(-m \lambda ^2 +2 m a^2 \sin^2 \theta)\sin^2
\theta}{A ^2 B ^{\frac{5}{2}} \lambda r}+O(r^{-3}).
 \end{aligned}
 \eeQ
Denote the frame $\breve{e_1}=\partial_r$,
$\breve{e_2}=\frac{\partial_\theta}{r}$,
$\breve{e_3}=\frac{\partial_\psi}{r\sin\theta}$ and $\{\breve{e}
^i\}$ the coframe. Using the asymptotic expansion $\rho ^2=r^2
+O(r)$, we find the components of the local angular momentum density
$\tilde h ^z $
 \beQ
 \begin{aligned}
\tilde h ^z (\breve{e}_2, \breve{e}_1)&=-\frac{3ma \sin \theta}{A^2 B^{\frac{5}{2}} r^2}+O(r^{-3}),\\
\tilde h ^z (\breve{e}_2, \breve{e}_3)&=\frac{m \lambda ^2 -2 m a^2
\sin^2 \theta}{A ^2 B ^{\frac{5}{2}} \lambda r^2}
+O(r^{-3}),\\
\tilde h ^z (\breve{e}_3, \breve{e}_2)&=\frac{m \lambda}{A ^2 B
^{\frac{3}{2}} r^2}+O(r^{-3}),\\
\tilde h ^z (\breve{e}_i, \breve{e}_j)&=O(r^{-3}).
 \end{aligned}
 \eeQ
Thus, in natural coordinates,
 \beQ
 \begin{aligned}
\tilde h ^z _{1r}&=\tilde h ^z (\breve{e}_1 \sin \theta \cos \psi
+\breve{e}_2 \cos \theta \cos \psi -\breve{e}_3 \sin \psi,
\breve{e}_1)\\
&=-\frac{3ma \sin \theta \cos \theta \cos \psi}{A^2 B^{\frac{5}{2}} r^2}+O(r^{-3}),\\
\tilde h ^z _{2r}&=\tilde h ^z (\breve{e}_1 \sin \theta \sin \psi
+\breve{e}_2 \cos \theta \sin \psi +\breve{e}_3 \sin \psi,\breve{e}_1)\\
&=-\frac{3ma \sin \theta \cos \theta \sin \psi}{A^2 B^{\frac{5}{2}} r^2}+O(r^{-3}),\\
\tilde h ^z _{3r}&=\tilde h ^z (\breve{e}_1 \cos \theta -\breve{e}_2
\sin \theta,\breve{e}_1)+O(r^{-3})\\
&=\frac{3ma \sin ^2 \theta}{A^2 B^{\frac{5}{2}} r^2}+O(r^{-3}).
 \end{aligned}
 \eeQ
Note that the range of $\bar \psi$ from $0$ to $2\pi$ gives that
 \beQ
-\frac{a}{\k ^2} \bar{t} \leq \psi \leq 2\big(1+\frac{a^2}{\k
^2}\big)\pi -\frac{a}{\k ^2} \bar{t},
 \eeQ
we obtain,
 \beQ
 \begin{aligned}
J_1 (z)&=\frac{A^2}{8\pi}\lim _{r \rightarrow \infty} \int
_{S_{r}}\tilde h ^z _{1r} *dr \\
&=-\frac{3ma}{8\pi}\int _{-\frac{a}{\k ^2} \bar{t}}
^{2\big(1+\frac{a^2}{\k ^2}\big)\pi -\frac{a}{\k ^2} \bar{t}} \int
_0 ^\pi
  \frac{\sin ^2 \theta \cos \theta \cos \psi}{B^{\frac{5}{2}}}d\theta d\psi =0,\\
J_2 (z)&=\frac{A^2}{8\pi}\lim _{r \rightarrow \infty} \int
_{S_{r}}\tilde h ^z _{2r} *dr\\
&=-\frac{3ma}{8\pi} \int _{-\frac{a}{\k ^2} \bar{t}}
^{2\big(1+\frac{a^2}{\k ^2}\big)\pi -\frac{a}{\k ^2} \bar{t}}
  \int _0 ^\pi \frac{\sin ^2 \theta \cos \theta \sin \psi}{B^{\frac{5}{2}}}d\theta d\psi =0,\\
J_3 (z)&=\frac{A^2}{8\pi}\lim _{r \rightarrow \infty} \int
_{S_{r}}\tilde h ^z _{3r} *dr\\
&=\frac{3ma}{8\pi} \int _{-\frac{a}{\k ^2} \bar{t}}
^{2\big(1+\frac{a^2}{\k ^2}\big)\pi -\frac{a}{\k ^2} \bar{t}}
  \int _0 ^\pi \frac{\sin ^3 \theta }{B^{\frac{5}{2}}}d\theta d\psi
  =\frac{ma}{1+\frac{a^2}{\lambda ^2}}.
 \end{aligned}
 \eeQ
\begin{rmk}
If we choose the range of $\psi$ from $0$ to $2\pi$, we obtain
 \beQ
J _1(z)=0,\ J_2(z)=0,\ J _3(z)=\frac{ma}{(1+\frac{a^2}{\lambda
^2})^2}.
 \eeQ
It is interesting that $J _3(z)$ conjugates to $J_{23}$ of \cite{HT}
by replacing $\k$ to $\sqrt{-1} \k$, i.e., the positive cosmological constant
to the negative cosmological constant.
\end{rmk}

\begin{rmk}
Note that the total energy and the total linear momentum vanish
for this $t$-slice. However, it does not contradict to the positive
mass theorem as it does not hold on the $t$-slice due to the
singularity. The situation is similar to the Schwarzschild spacetime
with negative mass. The positive mass theorem for black holes can
not apply to this spacetime as naked singularity occurs.
\end{rmk}

{\footnotesize {\it Acknowledgement.} The authors are indebted to
Ding Wang and Sibo Zheng for some useful conversations. M. Luo is
supported partially by the National Science Foundation of China
(10425525). N. Xie is supported partially by the National Science
Foundation of China (10801036). X. Zhang is supported partially by
the National Science Foundation of China (10421001, 10725105,
10731080), NKBRPC(2006CB805905) and the Innovation Project of
Chinese Academy of Sciences.}


\begin{thebibliography}{ADM}

\bibitem{AD} Abbott, L.F., Deser, S., Stability of gravity with a
cosmological constant, Nucl. Phys. B195(1982), 76-96.

\bibitem{ADM} Arnowitt, R., Deser, S., Misner, C.,
Coordinate invariance and energy expressions in general relativity,
Phys. Rev. 122(1961), 997-1006.

\bibitem{AMR} Astefanesei, D., Mann, R.B., Radu, E.,
Reissner-Nordstrom-de Sitter black hole, planar coordinates and
dS/CFT, J. High Ener. Phys. 0401(2004), 029.

\bibitem{BBM} Balasubramanian, V., de Boer, J., Minic, D., Mass, entropy,
and holography in asymptotically de Sitter spaces, Phys. Rev. D
65(2002), 123508.

\bibitem{B84} Bartnik, R., Existence of maximal surfaces in
asymptotically flat spacetimes, Commun. Math. Phys. 94(1984),
155-175.

\bibitem{B} Bartnik, R., The mass of an
asymptotically flat manifold, Comm. Pure. Appl. Math. 36(1986),
661-693.

\bibitem{BCM} Bartnik, R., Chru\'{s}ciel, P., \'{O} Murchadha, N.,
On maximal surfaces in asymptotically flat space-times, Commun.
Math. Phys. 130(1990), 95-109.

\bibitem{CK} Christodoulou, D., Klainerman, S.,
The global nonlinear stablity of Minkowski space, Princeton Math.
Series 41, Princeton Univ. Press, Princeton, 1993.

\bibitem{C1} Chru\'sciel, P., Boundary conditions at spatial
infinity from a Hamiltonian point of view, Topological Properties
and Global Structure of Space-Time (Erice, 1985), NATO, Adv. Sci.
Inst. Ser. B: Phys. 138, Plenum, New York, 1986, 49-59.

\bibitem{CH} Chru\'{s}ciel, P., Herzlich, M., The mass
of asymptotically hyperbolic Riemannian manifolds. Pacific J. Math.
212(2003), 231-264.

\bibitem{CGM} Clarkson, R., Ghezelbash, A.M., Mann, R.B., Mass,
action and entropy of Taub-Bolt-dS spacetimes, Phys. Rev. Lett.
91(2003), 061301.

\bibitem{Co} Corvino, J., Scalar curvature deformation and a gluing
construction for the Einstein constraint equations, Commun. Math.
Phys. 214(2000), 137--189.

\bibitem{GM} Ghezelbash, A.M., Mann, R.B., Entropy and mass bounds of
Kerr-de Sitter spacetimes, Phys. Rev. D 72(2005), 064024.

\bibitem{HE} Hawking, S., Ellis, G., The large scale structure of
space-time, Cambridge Univ. Press, Cambridge, 1973.

\bibitem{HT} Henneaux, M., Teitelboim, C., Asymptotically anti-de Sitter
spaces, Commun. Math. Phys. 98(1985), 391-424.

\bibitem{KT} Kastor, D., Traschen, J., A positive energy theorem for
asymptotically de Sitter spacetimes, Class. Quantum Gravity
19(2002), 5901-5920.

\bibitem{M} Maerten, D., Positive energy-momentum theorem for
AdS-asymptotically hyperbolic manifolds, Ann. Henri Poincar\'{e}
7(2006), 975-1011.

\bibitem{PT} Parker, T., Taubes, C., On
Witten's proof of the positive energy theorem, Commun. Math. Phys.
84(1982), 223-238.

\bibitem{RegT} Regge, T., Teitelboim, C., Role of surface
integrals in the Hamiltonian formulation of general relativity, Ann.
Phys. 88(1974), 286-318.

\bibitem{SY1} Schoen, R., Yau, S.T., On the proof of the
positive mass conjecture in general relativity, Commun. Math. Phys.
65(1979), 45-76.

\bibitem{SY2} Schoen, R., Yau, S.T., The energy and the
linear momentum of spacetimes in general relativity, Commun. Math.
Phys. 79(1981), 47-51.

\bibitem{SY3} Schoen, R., Yau, S.T., Proof of the positive
mass theorem II, Commun. Math. Phys. 79(1981), 231-260.

\bibitem{S} Shiromizu, T., Positivity of gravitational mass in
asymptotically de Sitter space-times, Phys. Rev. D 49(1994), 5026.

\bibitem{SIT} Shiromizu, T., Ida, D., Torii, T., Gravitational
energy, dS/CFT correspondence, and cosmic no-hair, J. High Ener.
Phys. 0010(2001), 010.

\bibitem{Wa}Wang, X., Mass for asymptotically hyperbolic manifolds.
J. Diff. Geom. 57(2001), 273-299.

\bibitem{W} Witten, E., A new proof of the positive energy
theorem, Commun. Math. Phys. 80(1981), 381-402.

\bibitem{W2} Witten, E., Quantum gravity in de Sitter space,
arXiv:hep-th/0106109.

\bibitem{XZ} Xie, N., Zhang, X., Positive mass theorems for
asymptotically AdS spacetimes with arbitrary cosmological constant,
Intern. J. Math. 19(2008), 285-302.

\bibitem{Z1} Zhang, X., Angular momentum and positive mass theorem,
Commun. Math. Phys. 206(1999), 137-155.

\bibitem{Z2} Zhang, X., Remarks on the total angular momentum in
general relativity, Commun. Theore. Phys., 39(2003), 521-524.

\bibitem{Z3}Zhang, X., A definition of total energy-momenta
and the positive mass theorem on asymptotically hyperbolic
3-manifolds I. Commun. Math. Phys. 249(2004), 529-548.
\end{thebibliography}
\end{document}